\documentclass[reqno,11pt]{amsart}

\usepackage{amsthm, amsmath, amssymb, stmaryrd}
\usepackage[utf8]{inputenc}
\usepackage[T1]{fontenc}

\usepackage[hyphenbreaks]{breakurl}
\usepackage[hyphens]{url}

\usepackage[shortlabels]{enumitem}
\usepackage[hidelinks]{hyperref}
\usepackage{microtype}

\usepackage{bm}
\usepackage[margin=1in]{geometry}

\usepackage[textsize=scriptsize,backgroundcolor=orange!5]{todonotes}

\usepackage[noabbrev,capitalize,sort]{cleveref}
\crefname{equation}{}{}
\numberwithin{equation}{section}

\newtheorem{theorem}{Theorem}[section]

\newtheorem{lemma}[theorem]{Lemma}

\theoremstyle{definition}

\theoremstyle{remark}

\newcommand\coin{
\mathrel{\ooalign{\hss$\bigcirc$\hss\cr\kern0.7ex\hbox{\scalebox{0.8}{$\$$}}}}\,}

\makeatletter
\newcommand\thankssymb[1]{\textsuperscript{\@fnsymbol{#1}}}
\makeatother

\author[Ting-Wei Chao]{Ting-Wei Chao\thankssymb{1}}
\author[Hung-Hsun Hans Yu]{Hung-Hsun Hans Yu\thankssymb{2}}

\thanks{\thankssymb{1}Department of Mathematics, Massachusetts Institute of Technology, Cambridge, MA, USA. Email: {\tt twchao@mit.edu}}

\thanks{\thankssymb{2}Department of Mathematics, Princeton University, Princeton, NJ 08544\@.  Email: {\tt hansonyu@princeton.edu}}

\title{A Purely Entropic Approach to the Rainbow Triangle Problem}

\begin{document}

\maketitle

\begin{abstract}
In this short note, we present a purely entropic proof that in a $3$-edge-colored simple graph with $R$ red edges, $G$ green edges, and $B$ blue edges, the number of rainbow triangles is at most $\sqrt{2RGB}$.
\end{abstract}

\section{Introduction}
Inspired by the multijoints problem, in \cite{CY23}, we considered the following graph-theoretic problem that is natural and interesting on its own.
Given a simple graph with $R$ red edges, $G$ green edges and $B$ blue edges, how many rainbow triangles can there be in the graph?
Using the entropy method, in that paper, we proved the following.
\begin{theorem}[\cite{CY23}]\label{thm:RT}
    Let $G=(V,E)$ be a simple graph, and each edge is colored in one of the three colors: red, green, or blue. Suppose there are $R$ red edges, $G$ green edges, $B$ blue edges, and $T$ rainbow triangles. Then we have 
    \[T^2\leq 2RGB.\]
\end{theorem}
We refer the readers back to the original paper \cite{CY23} for a more detailed account of the background and the connection to the multijoints problem.

Although the original proof uses the entropy method, what makes the proof different from other more well-known entropic proofs is that in the second part of the proof, entropy no longer plays a central role after all the terms are expressed in terms of probability of certain events.
In this note, we provide an alternative that stays in the realm of the entropy method throughout the entire proof.
Though it might be even more elusive how one could find this proof, it is arguably easier to digest conceptually as the key insight is now condensed into one single injection that is exhibited in \cref{subsection:keyinj}.

This note is structured as follows.
In \cref{section:RainbowTriangle}, we provide a self-contained entropic proof of the main theorem.
In \cref{section:conclusion}, we end this note with some remarks on the proof strategy and its possible implications.

\section{Main Proof}\label{section:RainbowTriangle}
In this section, we will prove \cref{thm:RT}.
Although the proof could fit into one page, to navigate the readers through various computational details, we split the proof up into three components so that it is easier to understand.

\subsection{Sampling random variables}
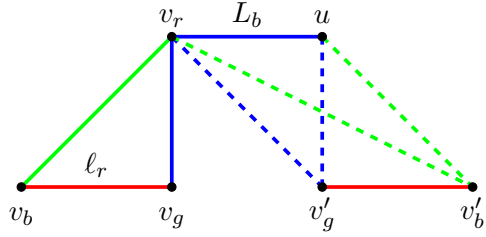
\begin{figure}[h]
    \centering
    
    \begin{tikzpicture}
    \coordinate (b) at (0,0);
    \coordinate (g) at (2,0);
    \coordinate (r) at (2,2);
    \coordinate (u) at (4,2);
    
    \coordinate (lr) at (1,0);
    
    \coordinate (lB) at (3,2);

    \coordinate (b') at (6,0);
    \coordinate (g') at (4,0);
    
    \draw[green, line width=1.4pt] (r) -- (b);
    \draw[red, line width=1.4pt] (b) -- (g);
    \draw[blue, line width=1.4pt] (g) -- (r);
    
    \draw[blue, line width=1.4pt] (r) -- (u);

    \draw[red, line width=1.4pt] (b') -- (g');

    \draw[blue, dashed, line width=1.4pt] (u) -- (g');
    \draw[green, dashed, line width=1.4pt] (u) -- (b');
    \draw[blue, dashed, line width=1.4pt] (r) -- (g');
    \draw[green, dashed, line width=1.4pt] (r) -- (b');
    
    \node[above] at (r) {$v_r$};
    \node[below,yshift=-4] at (g) {$v_g$};
    \node[below,yshift=-4] at (b) {$v_b$};
    
    \node[above] at (lr) {$\ell_r$};
    
    \node[above,yshift=1.5] at (u) {$u$};
    
    \node[above] at (lB) {$L_b$};
    
    \node[below] at (g') {$v'_g$};
    \node[below] at (b') {$v'_b$};

    \draw [fill] (r) circle (1.6pt);
    \draw [fill] (g) circle (1.6pt);
    \draw [fill] (b) circle (1.6pt);
    \draw [fill] (u) circle (1.6pt);
    \draw [fill] (g') circle (1.6pt);
    \draw [fill] (b') circle (1.6pt);
\end{tikzpicture}
    \caption{$v'_r$ is either $v_r$ or $u$.}
    \label{fig:rv}
\end{figure}

Throughout the proof, we will consider the random variables sampled as follows (see \cref{fig:rv}). 
Let $\Delta$ be a random rainbow triangle chosen uniformly at random, let $v_r,v_g,v_b$ be the vertices of $\Delta$ that are opposite to the red, green, blue edge in $\Delta$ respectively, and let $\ell_r=\{v_g,v_b\}$ be the red edge of $\Delta$. 
We sample a random vertex $u$ uniformly from the blue neighbors of $v_r$ such that $u$ and $\ell_r$ are conditionally independent given $v_r$, and we denote $L_b=\{v_r,u\}$.
Finally, consider the pair of random variables $(\Delta, L_b)$.
We resample a pair $(\Delta',L'_b)$ independently given $L'_b=L_b$ as a set.
Denote by $v'_r,v'_g,v'_b$ the vertices of $\Delta'$ opposite to the red, green, and blue edge in $\Delta'$ respectively. Note that we only require $L'_b=L_b$ as a set, so $v'_r$ could be either $v_r$ or $u$.

Note that $\Delta$ and $\Delta'$, when conditioning on $L_b$, are identically distributed and independent.
Therefore we have $H(\Delta \mid L_b) = H(\Delta'\mid L_b)$.
To upper bound $H(\Delta)$, it is natural to relate $H(\Delta)$ to $H(\Delta\mid L_b)$ and then replace $H(\Delta\mid L_b)$ with $H(\Delta'\mid L_b)$.
To do so, note that, by the chain rule, 
\[H(\Delta,u) = H(L_b)+H(\Delta\mid L_b)\]
and also
\[H(\Delta, u) = H(\Delta)+H(u\mid \Delta).\]
Notice that $H(u\mid \Delta) = H(u\mid v_r)$ as $u$ was conditionally independently drawn.
Therefore,
\[H(\Delta) = H(\Delta\mid L_b)+H(L_b)-H(u\mid v_r),\]
and hence
\[H(\Delta) = H(\Delta'\mid L_b)+H(L_b)-H(u\mid v_r).\]
As before, we would like to upper bound $2H(\Delta)$ in terms of entropies of the edges.
By the conditional independence of $\Delta$ and $\Delta'$ given $L_b$, it is clear that
\begin{align}\label{eq:keyentropy}
    2H(\Delta) = H(\Delta, \Delta'\mid L_b)+2H(L_b)-2H(u\mid v_r) = H(\Delta,u,\Delta')+H(L_b)-2H(u\mid v_r).
\end{align}
\subsection{Key injection}\label{subsection:keyinj}
Now we introduce the key component of the proof, which allows us to compress the information of $(v_b,v_g,u,\Delta')$ given $v_r$.
To do so, we will make several definitions.

We say that a tuple of three vertices $(x,y,z)$ \emph{forms a rainbow triangle in order} if $\{x,y\}$ is blue, $\{y,z\}$ is red, and $\{x,z\}$ is green. For any vertex $x$, let $N_g(x)$ and $N_b(x)$ be the set of green neighbours and blue neighbors of $x$ respectively. 
In addition, set
\[S_x=\left\{(x_1,x_2,y,z,z_1,z_2)\in V^6\middle|\begin{matrix}
    \text{$(x,x_1,x_2)$ and $(z,z_1,z_2)$ both form rainbow triangles in order,}\\
    \text{$\{x,y\}$ is blue, and $z\in\{x,y\}$},
\end{matrix}\right\},\]
and
\[T_x=\{(a,b_1,b_2,\ell)\in N_g(x)\times N_b(x)^2\times E\mid \text{$\ell$ is red}\}.\]
Notice that $S_x$ is defined so that $(v_b, v_g, u, v_r', v_g', v_b')$ is in $S_{v_r}$.
See below for an illustration of the definition of $T_x$.

\begin{figure}[h]
    \centering
    
    \begin{tikzpicture}
    \coordinate (b) at (0,0);
    \coordinate (g) at (2,0);
    \coordinate (r) at (2,2);
    \coordinate (u) at (4,2);
    
    \coordinate (l) at (5,0);

    \coordinate (b') at (6,0);
    \coordinate (g') at (4,0);
    
    \draw[green, line width=1.4pt] (r) -- (b);
    \draw[blue, line width=1.4pt] (g) -- (r);
    
    \draw[red, line width=1.4pt] (b') -- (g');
    \draw[blue, line width=1.4pt] (r) -- (u);

    \node[above,yshift=1.5] at (r) {$x$};
    \node[below] at (g) {$b_1$};
    \node[below,yshift=-2] at (b) {$a$};

    \node[above] at (u) {$b_2$};
    
    \node[above] at (l) {$\ell$};

    \draw [fill] (r) circle (1.6pt);
    \draw [fill] (g) circle (1.6pt);
    \draw [fill] (b) circle (1.6pt);
    \draw [fill] (u) circle (1.6pt);
    \draw [fill] (g') circle (1.6pt);
    \draw [fill] (b') circle (1.6pt);
\end{tikzpicture}
    \caption{The definition of $T_x$.}
    \label{fig:Tx}
\end{figure}
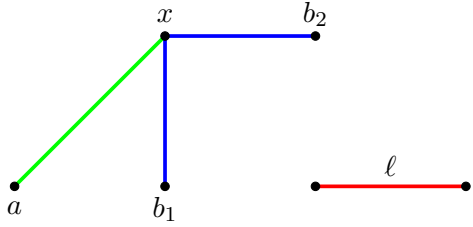

The key ingredient in this subsection is the following lemma.
\begin{lemma}\label{lemma:injection}
    For each $x\in V$, there is an injection $f_x:S_x\rightarrow T_x$.
\end{lemma}

We will prove the lemma later in the next subsection.
Here, we first see how the lemma implies the desired result.
\begin{proof}[Proof of \cref{thm:RT} assuming \cref{lemma:injection}]
    As $(v_g,v_b,u,v'_r,v'_g,v'_b)\in S_{v_r}$, we know that\\ $f_{v_r}(v_g,v_b,u,v'_r,v'_g,v'_b)$ is well-defined, which we denote as $(w,s_1,s_2,L_r)$.
    As $f_{v_r}$ is injective, it follows that 
    \begin{align}\label{eq:bijection}
        H(v_g,v_b,u,\Delta'\mid v_r)=H(w,s_1,s_2,L_r\mid v_r)
    \end{align}
    To combine it with \cref{eq:keyentropy}, we apply the chain rule to get that
    \begin{align*}
        2H(\Delta)= & H(v_b,v_g,u,\Delta'\mid v_r)+H(v_r)+H(L_b)-2H(u\mid v_r)\\
        =&H(w,s_1,s_2,L_r\mid v_r)+H(v_r)+H(L_b)-2H(u\mid v_r).
    \end{align*}

    By conditional subadditivity, the first term can be upper bounded by
    \[H(w\mid v_r)+H(s_1\mid v_r)+H(s_2\mid v_r)+H(L_r\mid v_r),\]
    and by the uniform bound, the choice of $u$ and monotonicty, we see that this can be upper bounded by
    \[H(w\mid v_r)+2H(u\mid v_r)+H(L_r)\]
    as $u$ is chosen uniformly from the blue neighbors of $v_r$.
    Hence,
    \[2H(\Delta)\leq H(w\mid v_r)+H(v_r)+H(L_r)+H(L_b) = H(v_r,w)+H(L_r)+H(L_b).\]
    As $\{v_r,w\}$ is a green edge, $L_r$ is a red edge and $L_b$ is a blue edge, by the uniform bound we get that
    \[\log_2 T^2 = 2H(\Delta)\leq \log_2(2G)+\log_2 R+\log_2 B = \log_2 2RGB,\]
    as desired.
\end{proof}

\subsection{Proof of the key injection}
Finally, in this subsection, we provide an elementary proof of \cref{lemma:injection}.

\begin{proof}[Proof of \cref{lemma:injection}]
    Fix a tuple $(x_1,x_2,y,z,z_1,z_2)\in S_x$. If $z=y$ and $(x,z_1,z_2)$ forms a rainbow triangle in order (see \cref{fig:case1}), then we define 
    \[f_x(x_1,x_2,y,z,z_1,z_2)=(z_2,z,z_1,\{x_1,x_2\}).\]

\begin{figure}[h]
    \centering
    
    \begin{tikzpicture}
    \coordinate (b) at (0,0);
    \coordinate (g) at (2,0);
    \coordinate (r) at (2,2);
    \coordinate (u) at (4,2);
    
    \coordinate (l) at (1,0);

    \coordinate (b') at (6,0);
    \coordinate (g') at (4,0);
    
    \draw[green, line width=1.4pt] (r) -- (b);
    \draw[red, line width=1.4pt] (b) -- (g);
    \draw[blue, line width=1.4pt] (g) -- (r);
    
    \draw[blue, line width=1.4pt] (r) -- (u);

    \draw[red, line width=1.4pt] (b') -- (g');

    \draw[blue, line width=1.4pt] (u) -- (g');
    \draw[green, line width=1.4pt] (u) -- (b');
    \draw[blue, line width=1.4pt] (r) -- (g');
    \draw[green, line width=1.4pt] (r) -- (b');
    
    \node[above,yshift=1.5] at (r) {$x$};
    \node[below,yshift=-4] at (g) {$x_1$};
    \node[below,yshift=-4] at (b) {$x_2$};

    \node[above] at (l) {$\ell$};
    
    \node[above] at (u) {$y=z=b_1$};

    \node[below] at (g') {$z_1=b_2$};
    \node[below,yshift=-2] at (b') {$z_2=a$};

    \draw [fill] (r) circle (1.6pt);
    \draw [fill] (g) circle (1.6pt);
    \draw [fill] (b) circle (1.6pt);
    \draw [fill] (u) circle (1.6pt);
    \draw [fill] (g') circle (1.6pt);
    \draw [fill] (b') circle (1.6pt);
\end{tikzpicture}
    \caption{The case when $z=y$ and $(x,z_1,z_2)$ forms a rainbow triangle in order.}
    \label{fig:case1}
\end{figure}
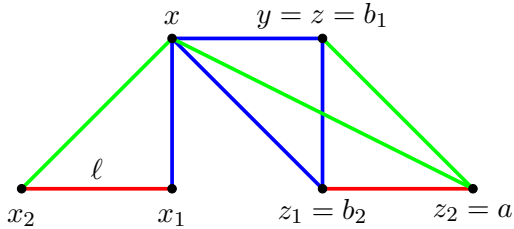

If $z=y$ and $(x,z_2,z_1)$ forms a rainbow triangle in order (see \cref{fig:case2}), then we define 
\[f_x(x_1,x_2,y,z,z_1,z_2)=(z_1,z,z_2,\{x_1,x_2\}).\]

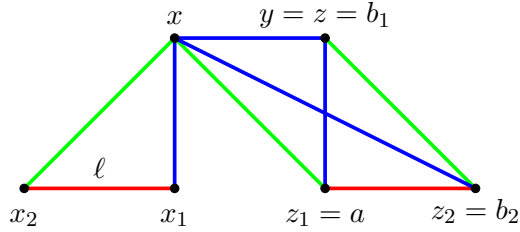
\begin{figure}[h]
    \centering
    
    \begin{tikzpicture}
    \coordinate (b) at (0,0);
    \coordinate (g) at (2,0);
    \coordinate (r) at (2,2);
    \coordinate (u) at (4,2);
    
    \coordinate (l) at (1,0);

    \coordinate (b') at (6,0);
    \coordinate (g') at (4,0);
    
    \draw[green, line width=1.4pt] (r) -- (b);
    \draw[red, line width=1.4pt] (b) -- (g);
    \draw[blue, line width=1.4pt] (g) -- (r);
    
    \draw[blue, line width=1.4pt] (r) -- (u);

    \draw[red, line width=1.4pt] (b') -- (g');

    \draw[blue, line width=1.4pt] (u) -- (g');
    \draw[green, line width=1.4pt] (u) -- (b');
    \draw[green, line width=1.4pt] (r) -- (g');
    \draw[blue, line width=1.4pt] (r) -- (b');
    
    \node[above,yshift=1.5] at (r) {$x$};
    \node[below,yshift=-4] at (g) {$x_1$};
    \node[below,yshift=-4] at (b) {$x_2$};

    \node[above] at (l) {$\ell$};
    
    \node[above] at (u) {$y=z=b_1$};

    \node[below,yshift=-4] at (g') {$z_1=a$};
    \node[below] at (b') {$z_2=b_2$};

    \draw [fill] (r) circle (1.6pt);
    \draw [fill] (g) circle (1.6pt);
    \draw [fill] (b) circle (1.6pt);
    \draw [fill] (u) circle (1.6pt);
    \draw [fill] (g') circle (1.6pt);
    \draw [fill] (b') circle (1.6pt);
\end{tikzpicture}
    \caption{The case when $z=y$ and $(x,z_2,z_1)$ forms a rainbow triangle in order.}
    \label{fig:case2}
\end{figure}

Otherwise (see \cref{fig:case3}), we define
\[f_x(x_1,x_2,y,z,z_1,z_2)=(x_2,x_1,y,\{z_1,z_2\}).\]

\begin{figure}[h]
    \centering
    
    \begin{tikzpicture}
    \coordinate (b) at (0,0);
    \coordinate (g) at (2,0);
    \coordinate (r) at (2,2);
    \coordinate (u) at (4,2);
    
    \coordinate (l) at (5,0);

    \coordinate (b') at (6,0);
    \coordinate (g') at (4,0);
    
    \draw[green, line width=1.4pt] (r) -- (b);
    \draw[red, line width=1.4pt] (b) -- (g);
    \draw[blue, line width=1.4pt] (g) -- (r);
    
    \draw[blue, line width=1.4pt] (r) -- (u);

    \draw[red, line width=1.4pt] (b') -- (g');

    \node[above,yshift=1.5] at (r) {$x$};
    \node[below] at (g) {$x_1=b_1$};
    \node[below,yshift=-3] at (b) {$x_2=a$};

    \node[above] at (l) {$\ell$};
    
    \node[above] at (u) {$y=b_2$};

    \node[below,yshift=-3] at (g') {$z_1$};
    \node[below,yshift=-3] at (b') {$z_2$};

    \draw [fill] (r) circle (1.6pt);
    \draw [fill] (g) circle (1.6pt);
    \draw [fill] (b) circle (1.6pt);
    \draw [fill] (u) circle (1.6pt);
    \draw [fill] (g') circle (1.6pt);
    \draw [fill] (b') circle (1.6pt);
\end{tikzpicture}
    \caption{The remaining case.}
    \label{fig:case3}
\end{figure}
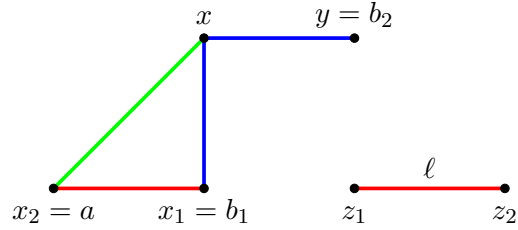

It is not hard to check that in all three cases $f_x(x_1,x_2,y,z,z_1,z_2)\in T_x$. It remains to show that $f_x$ is injective. Given any $(a,b_1,b_2,\ell)\in T_x$ in the image, observe that the color of the edge $\{a,b_1\}$ indicates the case which $(a,b_1,b_2,\ell)$ is mapped from. Namely, if $\{a,b_1\}$ is green (blue, red, resp.) then it must come from the first (second, third, resp.) case above. Therefore, when $\{a,b_1\}$ is green, the only tuple that can be in the preimage is 
\[(x_1,x_2,y,z,z_1,z_2)=(c_1,c_2,b_1,b_1,b_2,a),\]
where $\ell=\{c_1,c_2\}$ and the order of $c_1,c_2$ is chosen such that $(x,c_1,c_2)$ forms a rainbow triangle in order. Similarly, when $\{a,b_1\}$ is blue, the only tuple that can be in the preimage is 
\[(x_1,x_2,y,z,z_1,z_2)=(c_1,c_2,b_1,b_1,a,b_2).\]
Finally, when $\{a,b_1\}$ is red, we know that $(x_1,x_2,y)=(b_1,a,b_2)$. If $x$ and $\ell$ forms a rainbow triangle, then we must have $z=x$, as otherwise it must have fallen into the first two cases instead. 
If $x$ and $\ell$ do not form a rainbow triangle, then we must have $z=y$ since $z$ and $\ell$ forms a rainbow triangle in this case. 
Lastly, $z_1,z_2$ can be deduced uniquely as well as they must be the two vertices of $\ell$ such that $(z,z_1,z_2)$ forms a rainbow triangle in order.
\end{proof}

\section{A concluding remark}\label{section:conclusion}
In the proof, the most important ingredient seems to be the key injection (\cref{lemma:injection}).
It seems to us that this injection corresponds to the case analysis of type A and type B triangles in the original proof \cite{CY23}.
It is interesting to think about how one might discover the injection and find use of it without knowing the previous proof.
One possible thought is that, perhaps, there can be programs that brute-force through useful configurations while keeping track of the relations between their ``entropies'', similar to the philosophy behind flag algebras.
This may not be such a surprise as a related problem was previously studied via flag-algebraic computations \cite{BHLPVY17}.
There, the number of vertices is given instead of the number of edges of each color.
The idea of emulating flag algebra using the entropy method was also communicated privately by Yufei Zhao to us.
Although we do not see at this point how to put this thought into practice, we think it is an interesting idea that is worth more attention both because of the various generalizations of the rainbow triangle problems, and because of the many open problems that might be related to the entropy method.

\section*{Acknowledgement}
The proof was discovered when the authors were hosted by Yufei Zhao at MIT.
The authors would like to thank him for his hospitality and the inspiring discussions that lead to the proof in this note.

\bibliographystyle{amsplain0}
\bibliography{ref_KK}
\end{document}